\documentclass[11pt]{article}
\usepackage{amssymb,amsfonts,amsmath,latexsym,epsf,tikz,url}
\newtheorem{theorem}{Theorem}[section]

\newtheorem{lemma}[theorem]{Lemma}

\newcommand{\proof}{\noindent{\bf Proof.\ }}
\newcommand{\qed}{\hfill $\square$\medskip}

\begin{document}

\title{On the roots of total domination polynomial of graphs, II}

\author{Saeid Alikhani$^{}$\footnote{Corresponding author} and Nasrin Jafari }

\date{\today}

\maketitle

\begin{center}

   Department of Mathematics, Yazd University, 89195-741, Yazd, Iran\\
{\tt alikhani@yazd.ac.ir}\\

\end{center}


\begin{abstract}
Let $G = (V, E)$ be a simple graph of order $n$. The total dominating set of $G$ is a subset $D$ of $V$ that every vertex of $V$ is adjacent to some vertices of $D$. The total domination number of $G$ is equal to minimum cardinality of  total dominating set in $G$ and is denoted by $\gamma_t(G)$. The total domination polynomial of $G$ is the polynomial $D_t(G,x)=\sum_{i=\gamma_t(G)}^n d_t(G,i)x^i$, where $d_t(G,i)$ is the number of total dominating sets of $G$ of size $i$. A root of $D_t(G, x)$ is called a total  domination root of $G$. The set of total domination roots of graph $G$ is denoted by $Z(D_t(G,x))$. In this paper we show that $D_t(G,x)$ has $\delta-2$ non-real roots and if all roots of $D_t(G,x)$ are real then $\delta\leq 2$, where $\delta$ is the minimum degree of vertices of $G$. Also we show that if $\delta\geq 3$ and $D_t(G,x)$ has exactly three distinct roots, then
$Z(D_t(G,x))\subseteq \{0, -2\pm \sqrt{2}i, \frac{-3\pm \sqrt{3}i}{2}\}$.  Finally we study the location roots of total domination polynomial of some families of graphs.  
\end{abstract}

\noindent{\bf Keywords:} total domination polynomial, root, real.

\medskip
\noindent{\bf AMS Subj.\ Class.:}  05C12.

\section{Introduction}
Let $G = (V, E)$ be a simple graph. The order of $G$ is the number of vertices of $G$. For any vertex $ v \in V$, the open neighborhood of $v$ is the set $N(v)=\{ u \in V | uv \in E\}$ and the closed neighborhood is the set $N[v]=N (v) \cup \{v\}$.
For a set $S\subset V$, the open neighborhood of $S$ is the set $N(S)=\bigcup_{v\in S }N(v)$ and the closed neighborhood of $S$ is the set $N[S]=N (S) \cup S$. 
A leaf (end-vertex) of a graph is a vertex of degree one, while a support vertex is a vertex adjacent to a leaf. 
The set $D\subset V$ is a total dominating set if every vertex of $V$ is adjacent to some vertices of $D$, or equivalently, $N(D)=V$. The total domination  number $\gamma_t(G)$ is the minimum cardinality of a total dominating set in $G$. A total dominating set with cardinality $\gamma_t(G)$ is called a $\gamma_t$-set. An $i$-subset of $V$ is a subset of $V$ of cardinality $i$. Let $\mathcal{D}_t(G, i)$ be the family of total dominating sets of $G$ which are $i$-subsets and let $d_t(G,i)=|\mathcal{D}_t(G, i)|$. The polynomial $D_t(G; x)=\sum_{i=1}^n d_t(G,i)x^i$ is defined as total domination polynomial of $G$.
As an example,  $D_t(K_n,x) = (x + 1)^n-nx-1$ and $D_t(K_{1,n},x)=x((x+1)^n-1)$. 
A root of $D_t(G, x)$ is called a total  domination root of $G$. The set of total domination roots of graph $G$ is denoted by $Z(D_t(G,x))$. 
For many graph polynomials, their roots have attracted considerable attention. 
For example in \cite{BHN}  Brown, Hickman, and Nowakowski proved that the real roots of the independence polynomials are dense in the interval $(-\infty,0]$, while the complex roots are dense in the complex plane. For matching polynomial, in \cite{Heil}   was proved that all roots of the matching polynomials are real. Also it was shown that if a graph has a Hamiltonian path, then all roots of its matching polynomial are simple (see Theorem 4.5 of [15]).  For domination polynomial,  Brown and Tufts in \cite{Brown} studied the location of domination roots and they proved that the set of all domination roots is dense in the complex plane. For graphs with few domination roots see \cite{euro}. 
Related to the roots of total domination polynomials there are a few papers. See \cite{nasrin,nasrin3} for more details. Recently authors  in \cite{nasrin} shown that all roots of $D_t(G, x)$ lie in the circle with center $(-1, 0)$ and 
radius $\sqrt[\delta]{2^n-1}$, where $\delta$ is the minimum degree of $G$ and $n$ is the order of $G$. As a consequence, they proved that if $\delta\geq \frac{2n}{3}$,
then every integer root of $D_t(G, x)$ lies in the set $\{-3,-2,-1,0\}$.

\medskip 
In this paper  we show that $D_t(G,x)$ has $\delta-2$ non-real roots and if all roots of $D_t(G,x)$ are real then $\delta\leq 2$. Also we show that if $\delta\geq 3$ and $D_t(G,x)$ has exactly three distinct roots, then
$Z(D_t(G,x))\subseteq \{0, -2\pm \sqrt{2}i, \frac{-3\pm \sqrt{3}i}{2}\}$.
Finally we study the location roots of total domination polynomial of some families of graphs.

\section{Main results}

\vglue-10pt
\indent

In this section we obtain more results on total domination roots. 
Oboudi in \cite{Oboudi} has studied graphs whose domination polynomials have only real roots. More precisely he obtained the number of non-real roots of domination polynomial of graphs. Similarly, we do it for total domination roots, in the next theorem. 

\begin{theorem}
	Let G be a connected graph of order $n\geq 2$. 
	\begin{enumerate}
		\item[i)] If all roots of $G$ are real, then $\delta=1 ~~or~~2$.
		\item[ii)]The polynomial $D_t(G,x)$ has at least $\delta-2$ non-real roots.
	\end{enumerate}\label{7}
\end{theorem}

\proof
Let $g(x) = D_t(G, x)$ and $g^{(m)}(x)$ be the $m$-th derivative of $g(x)$ with respect to $x$. It is easy to see that if $i \geq n-\delta+1$, then $d_t(G,i)={n \choose i}$ and if $i\leq n-\delta$, then $d_t(G,i)<{n \choose i}$, where $d_t(G,i)$ is the number of total dominating sets of $G$ with cardinality $i$, for every natural number $i$. Thus there exists a polynomial $f (x)$ with positive coefficients and with degree $n-\delta$ such that $D_t(G, x) = (x+1)^n-f (x)$. Since all roots of $g(x)$ are real, by Rolle's theorem we conclude that all roots of $g^{(n-\delta)}(x)$ are real as well. On the other hand $g^{(n-\delta)}(x) = \frac{n!} {\delta!} (x + 1)^{\delta} - a(n - \delta)!$, where $a$ is the coefficient of $x^{n-\delta}$ in $f (x)$. Since all roots of $g^{(n-\delta)}(x)$ are real, this shows that $\delta \leq 2$. Since $G$ is connected, so $\delta = 1$ or $2$.

Now suppose that $g(x)$ has exactly $r$ real roots. Using Rolle's theorem one can see that $g^{(n-\delta)}(x)$ has at least $r-(n-\delta)$ real roots. On the other hand $g^{(n-\delta)}(x) = \frac{n!}{\delta!} (x + 1)^{\delta} - a(n-\delta)!$. 
Thus $r-(n-\delta)\leq 2$. Therefore $g(x)$ has at least $\delta-2$ non-real roots.\qed

\begin{theorem}{\rm{\cite{nasrin3}}}\label{8}
	If $G=(V,E)$ is a graph of order $n$ with $r$ support vertices, then $d_t(G,n-1)=n-r$.
\end{theorem}

\begin{theorem}{\rm \cite{Hening}}\label{gamma_t}
	If $G$ is a graph of order $n$ with $\delta(G) \geq 3$, then $\gamma_t(G)\leq \frac{n}{2}$.
\end{theorem}

The study of graphs which their polynomials have few roots can give sometimes a surprising information about the structure of the graph. If $A$ is the adjacency matrix of $G$, then the eigenvalues of $A$, $\lambda_1 \geq \lambda_2\geq \ldots \geq \lambda_n$ are said to be the eigenvalues of the graph $G$. These are the roots of the characteristic polynomial $\phi(G,\lambda)=\prod_{i=1}^n (\lambda-\lambda_i).$ For more details on the characteristic polynomials. The characterization of graphs with few distinct roots of characteristic polynomials (i.e. graphs with few distinct eigenvalues) have been the subject of many researches. Graphs with three adjacency eigenvalues have been studied by Bridges and Mena \cite{bridges} and Klin and Muzychuk \cite{muz}. Also van Dam studied graphs with three and four distinct eigenvalues \cite{dam1,dam2,dam3,dam4}. Graphs with three distinct eigenvalues and index less than $8$ were studied by Chuang and Omidi in \cite{omidi}. Graphs with few domination roots were studied in \cite{euro}.  In \cite{nasrin3}, authors studied graphs with exactly two total domination roots $\{−3,0\}$, $\{−2,0\}$ and $\{−1,0\}$. Here we study graphs with three 
distinct total domination roots.

\begin{theorem}
	Let $G$ be a graph with $\delta\geq 3$. If $D_t(G,x)$ has exactly three distinct roots, then
	\[
	Z(D_t(G,x))\subseteq \{0, -2\pm \sqrt{2}i, \frac{-3\pm \sqrt{3}i}{2}\}
	\]
\end{theorem}

\proof
Let $G$ be a connected graph of order $n$ and $Z(D_t(G,x))=\{0,a,b\}$ that $a\neq b$. Therefore $D_t(G,x)=x^i(x-a)^j(x-b)^k$, 
for some $i, j, k$. So by Theorem \ref{8}, we have

\begin{equation}
	-(ja+kb)=n.\label{e1}
\end{equation} 

Also because $d_t(G, i)={n \choose i}$ for $i\geq n-\delta+1$, we have 

\begin{equation}
	{j \choose 2}a^2+{k \choose 2}b^2+jkab=d_t(G,n-2)={n \choose 2}.\label{e2}
\end{equation}

Let $P(x)$ be the minimal polynomial of $a$ over $\mathbb{Q}$. Clearly, all roots of $P(x)$ are simple. This implies that $deg(P(x))=1~ or ~2$. We consider two cases.
\begin{enumerate}
	\item[Case 1.] $deg(P(x))=1$. So $D_t(G,x)=x^i(x-a)^j(x-b)^k$, where $-a, -b\in \mathbb{N}$. By Theorem \ref{7}, we have $\delta=1~or~2$, a contradiction.
	
	\item[Case 2.] $deg(P(x))=2$. In this case since $D_t(G,x)$ has three distinct roots, the minimal polynomial of $b$ over $\mathbb{Q}$ is also $P(x)$, Thus we have $D_t(G,x)=x^i(x^2+rx+s)^j$, where $P(x)=x^2+rx+s$. We have $i+2j=n$, and also by (\ref{e1}), $-(a+b)j=n$. By Theorem \ref{gamma_t}, $i\leq \frac{n}{2}$. Therefore $j\geq \frac{n}{4}$. Since $-(a+b)j=n$ and $a+b$ is an integer, we have $-(a+b)\in \{1, 2, 3, 4 \}$. We consider four cases:
	\begin{enumerate}
		\item[Subcase 2.1.] If $a+b=-1$, then $j=n$, a contradiction.
		
		\item[Subcase 2.2.] If $a+b=-2$, then $j=\frac{n}{2}$, a contradiction.
		
		\item[Subcase 2.3.] If $a+b=-3$, then $i=j=\frac{n}{3}$, so we have $D_t(G,x)=x^{\frac{n}{3}}(x^2+rx+s)^{\frac{n}{3}}$. Now, by (\ref{e2}) we have
		
		\[{\frac{n}{3} \choose 2}(a^2+b^2)+\frac{n^2ab}{9}={n \choose 2}.\]
		
		In the other hand, since $a+b=-3$, we conclude that $a^2+b^2=9-2ab$. Thus by simple calculation we obtain $nab=3n$. Therefore $ab=3$. By using $a+b=-3$, we have 
		\[
		a\in \{\frac{-3\pm \sqrt{3}i}{2}\}
		\]
		
		\item[Subcase 2.4.] Now, suppose that $a+b=-4$. Then $i=\frac{n}{2}$ and $j=\frac{n}{4}$. With the same calculations, we have $ab=6$. Using the fact that $a+b=-4$, we have
		$a\in \{ -2\pm \sqrt{2}i\}.$\qed
	\end{enumerate}
\end{enumerate}

As noted before, in \cite{nasrin3}, authors studied graphs with exactly two total domination roots $\{-3,0\}$, $\{-2,0\}$ and $\{-1,0\}$. Here we present 
a family of graphs whose total domination roots are $-1$ and $0$. 

\begin{figure}[htb]
	\centering
	\includegraphics[height=3.4cm]{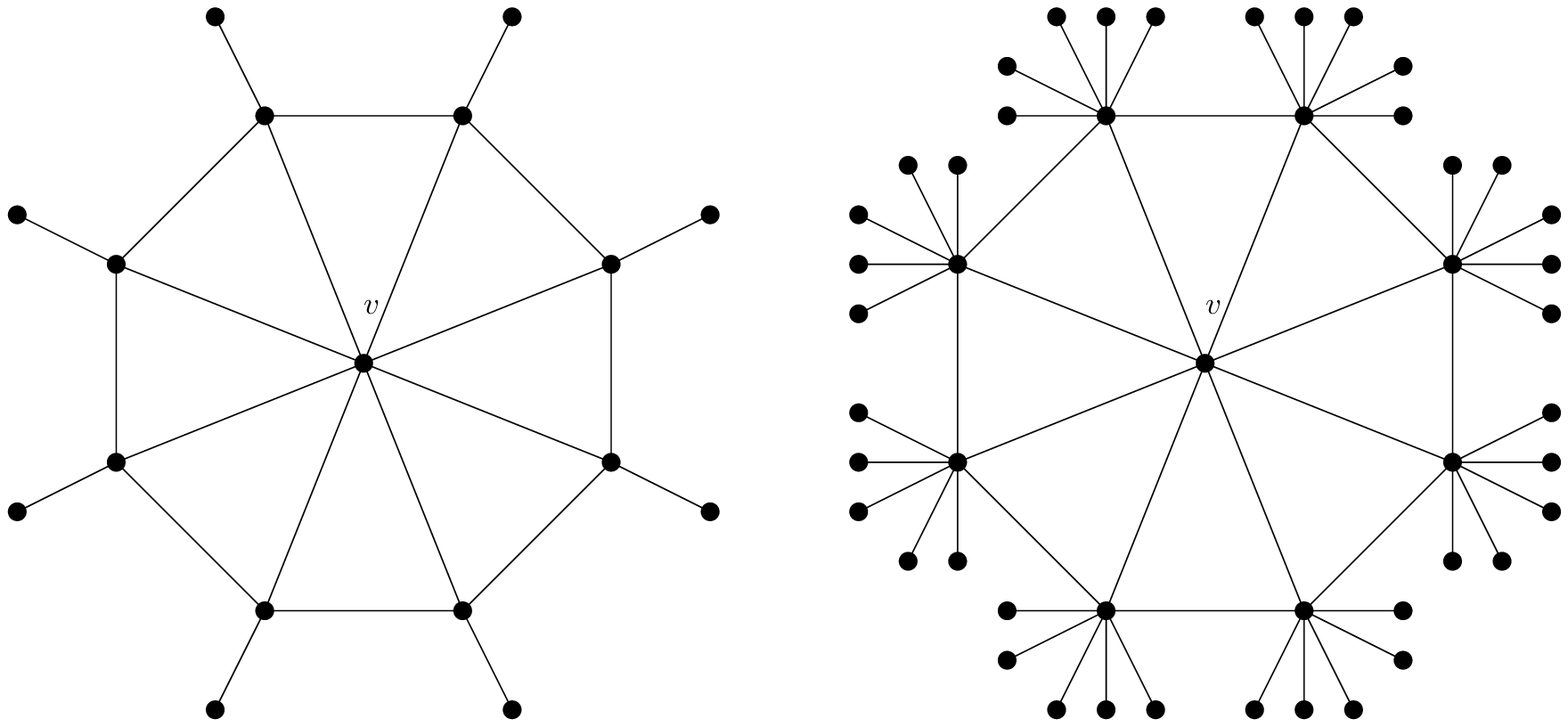}
	\caption{Helm graph $H_8$ and generalized helm graph $H_{8,5}$, respectively.}\label{Helm}
\end{figure}

The helm graph $H_n$ is obtained from the wheel graph $W_n$ by attaching a pendent edge at each vertex of the $n$-cycle of the wheel. We define generalized helm graph $H_{n,m}$, the graph is obtained from the wheel graph $W_n$ by attaching $m$ pendent edges at each vertex of the $n$-cycle of the wheel.
We recall that corona product of two graphs $G$ and $H$ is denoted by $G\circ H$ 
and was introduced by Harary \cite{Harary, Harar}. This graph formed from one copy of $G$ and $|V(G)|$ copies of $H$, where the $i$-th vertex of $G$ is adjacent to every vertex in the $i$-th copy of $H$. 
We need the following theorems:

\begin{theorem}{\rm\cite{Dod}}\label{10}
	Let $G = (V,E)$ be a graph and $u,v\in V$ two non-adjacent vertices of the graph with $N(u)\subseteq N(v)$. Then
	\[ D_t(G, x) = D_t(G \setminus v, x) + xD_t(G/v, x) +x^2\sum_{w\in N(v)\cap N(u)} D_t(G \setminus N[\{v,w\}], x).\]
\end{theorem}

\begin{theorem}{\rm\cite{nasrin}}\label{corona}
	For any graph $G$ of order $n\geq 2$, $D_t(G\circ \overline{K_m}, x)=x^n(1+x)^{mn}$.
\end{theorem}

\begin{theorem}
	For every natural number $n, m$, we have
	\begin{enumerate}
		\item[i)]$D_t(H_n,x)=x^n(x+1)^{n+1}$,
		\item[ii)]$D_t(H_{n,m},x)=x^n(1+x)^{mn+1}$.
	\end{enumerate} 
\end{theorem}
\proof
Let $v$ be the center vertex of wheel in helm graph $H_n$ and $H_{n,m}$. By Theorems \ref{10} and  \ref{corona} we have 

\begin{enumerate}
	\item[i)]$D_t(H_n,x)=D_t(C_n\circ K_1,x)+xD_t(K_n\circ K_1,x)=(1+x)(x(1+x))^n$,
	\item[ii)]$D_t(H_{n,m},x)=D_t(C_n\circ \overline{K_m},x)+xD_t(K_n\circ \overline{K_m},x)=(1+x)(x(1+x)^m)^n$.
\end{enumerate}
So we have the result.\qed

The lexicographic product is also known as graph substitution, a name that bears witness to the fact that $G[H]$ can be obtained from $G$ by substituting a copy $H_u$ of $H$ for every vertex $u$ of $G$ and then joining all vertices of $H_u$ with all vertices of $H_v$ if $\{u, v\} \in E(G)$.

\begin{theorem}
	Let $K_m$, $K_n$ be complete graphs of order $m$ and $n$. The total domination polynomial of lexicographic product of $K_m$ and $K_n$ is 
	\[D_t(K_m[K_n],x)=D_t(K_m,D(K_n,x))+mD_t(K_n,x)\]
\end{theorem}
\proof
Note that $K_m[K_n]\cong K_{mn}$, So the result is obtained.\qed

\begin{figure}[htb]
	\centering
	\includegraphics[height=5.25cm]{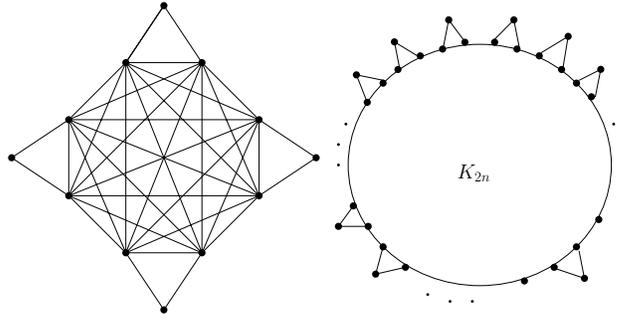}
	\caption{Graphs $G_4$ and $G_n$ in proof of theorem \ref{thFn,4}, respectively.}
	\label{Gn}
\end{figure}

\begin{figure}[htb]
	\centering
	\includegraphics[height=3cm]{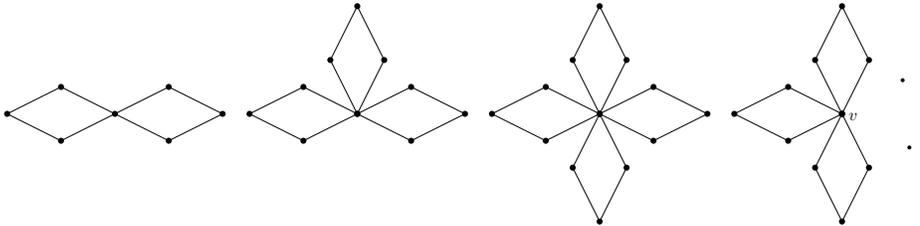}
	\caption{Friendship graphs $F_{2,4}, F_{3,4}, F_{4,4}$ and $F_{n,4}$, respectively.}
	\label{Dutch}
\end{figure}

The generalized friendship graph $F_{n,q}$ is a collection of $n$ cycles (all of order $q$), meeting at a common vertex (see Figure \ref{Dutch}). The generalized friendship graph may also be referred to as a flower \cite{Schiermeyer}. For $q=3$ the graph $F_{n,q}$ is denoted simply by $F_n$ and is friendship graph. The total domination polynomial of $F_n$ and its roots studied in \cite{nasrin}. Here, we compute the total domination number of $F_{n,4}$. 
To study the total domination roots of $F_{n,4}$ we first prove the following theorem which is a formula for the total domination polynomial of graph $G_n$ depicted in Figure \ref{Gn}.
\begin{theorem}\label{new}
	For any $n\in \mathbb{N}$, $D_t(G_n,x)=(x(x+1)(x+2))^n$.
\end{theorem} 
\proof
It is enough to prove the following equation for all $r\geq n$:
\begin{center}
	$d_t(G_n,r)=\sum \limits_{k=0}^{r-n}{n \choose r-n-k}{n \choose k}2^{n-k}$.
\end{center}
Suppose that $\{v_1,v_2, \ldots, v_{2n}\}$ be the set of vertecis of $K_{2n}$ and $\{u_1, u_2, \ldots, u_n\}$ are other vertecis of $G_n$. Now let $V_1=\{v_1, v_3, \ldots, v_{2n-1}\}$, $V_2=\{v_2,v_4, \ldots, v_{2n}\}$ and $U=\{u_1, u_2, \ldots, u_n\}$. For every $r\geq n$ and every $0\leq k\leq r-n$, choose $k$ vertecis of $V_1$ by ${ n \choose k}$ manner, therefore $k$ vertices of $U$ are total dominated. Now for total dominating $n-k$ vertecis of $U$, we can choose  $2^{n-k}$ vertices  of $V_1$ and $V_2$ and finally choose $r-n-k$ vertecis of $U$ by ${n \choose r-n-k}$ manner. So we have the result.\qed

\begin{theorem}\label{thFn,4}
	For every natural number $n$, total domination polynomial of generalize friendship graph $F_{n,4}$ is
	$$D_t(F_{n,4},x)=x^{n+1}(x+2)^n((x+1)^n+x^{n-1}).$$
\end{theorem}
\proof
Let $v$ be center vertex of $F_{n,4}$. By theorem \ref{10} we have
$$D_t(F_{n,4},x)=(D_t(P_3,x))^n+xD_t(G_n,x)$$
where $G_n$ is graph in Figure \ref{Gn} and so by Theorem \ref{new} we have the result.\qed

We need the following lemma to obtain more results:
\begin{lemma}{\rm\cite{Brown}} \label{limit}
	$ lim_{n \rightarrow \infty} ln(n)\Big(\frac{ln(n)-1}{{ ln(n)}}\Big)^n=0.$
\end{lemma}

The basic idea of the following result follows from the proof of Theorem $8$ in \cite{Brown}. 

\begin{theorem} 
	\begin{enumerate}
		\item[i)] For every natural even number $n$, the total domination polynomial of the generalized friendship graph, $D_t(F_{n,4}, x)$, where $n\geq2$, has a real root in the interval $(-1,0)$
		\item[ii)] The total domination polynomial of the generalized friendship graph, $D_t(F_{n,4}, x)$, 
		where $n\geq2$,  has a real root in the interval $(-n,-ln(n))$, for $n$ sufficiently large.
	\end{enumerate}
\end{theorem} 
\proof
\begin{enumerate}
	\item[i)]Let $f(x)=(x+1)^n+x^{n-1}$. So $f(0)=1$ and $f(-1)=(-1)^{n-1}=-1$. By the intermediate value theorem, we have result.
	\item[ii)] Suppose that $$f_{2n}(x)=x^{n+1}((x+1)^n+x^{n-1}).$$
	Observe that 
	$$f_{2n}(x)=x^{2n+1}+(n+1)x^{2n}+{n\choose n-2}x^{2n-1}+{n\choose n-3}x^{2n-2}+...+nx^{n+2}+x^{n+1}.$$
	Consider 
	
	$$f_{2n}(-n)=(-1)^{2n+1}n^{2n+1}\Big(1-\frac{n+1}{n}+\frac{{n\choose 2}}{(n)^2}-...+\frac{(-1)^n}{(n)^{{n}}}\Big).$$ 
	
	So $f_{2n}(-n)<0$ for $n$ sufficiently large, because the following inequality is true for $n$ sufficiently large, 
	
	$$ \frac{n+1}{n}-\frac{{n\choose 2}}{(n)^2}+...-\frac{(-1)^n}{(n)^{{n}}}<1.$$ 
	
	Now consider
	\begin{eqnarray*}
		f_{2n}(-ln(n)) &=& (-ln(n))^{n+1}(1-ln(n))^{n}+(-ln(n))^{2n}\\ 
		&=&(ln(n))^{2n}\Big(1-ln(n)\Big(\frac{ln(n)-1}{ln(n)}\Big)^{n}\Big).
	\end{eqnarray*}
	
	From Lemma \ref{limit}, we have $ ln(n)\Big(\frac{ln(n)-1}{{ln(n)}}\Big)^n \rightarrow 0$, as $n\rightarrow \infty$ which implies that $f_{2n}(-ln(n))>0$. By the Intermediate Value Theorem, for sufficiently large
	$n$, $f_{2n}(x) = D_t(F_n, x)$ has a real root in the interval $(-n,-ln(n))$.\qed
\end{enumerate}

\medskip

Figure \ref{Fn,4} shows the total domination roots of $F_{n,4}$ for $2 \leq n \leq 30$. 

\begin{figure}[htb]
	\centering
	\includegraphics[height=8cm]{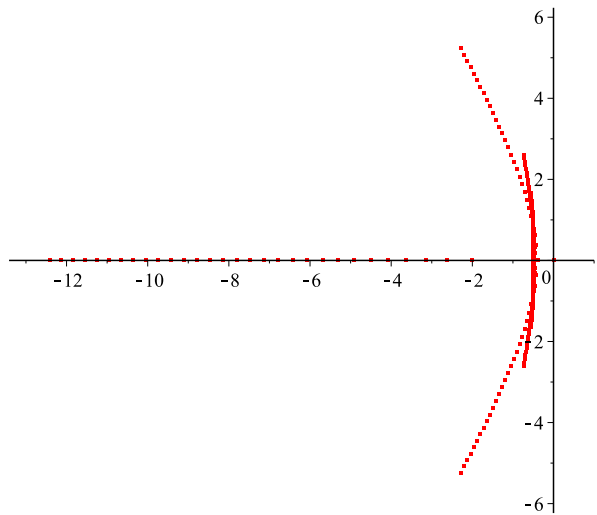}
	\caption{Total domination roots of $F_{n,4}$, for $2 \leq n \leq 30$.}\label{Fn,4}
\end{figure}

\begin{figure}[htb]
	\centering
	\includegraphics[height=6.4cm]{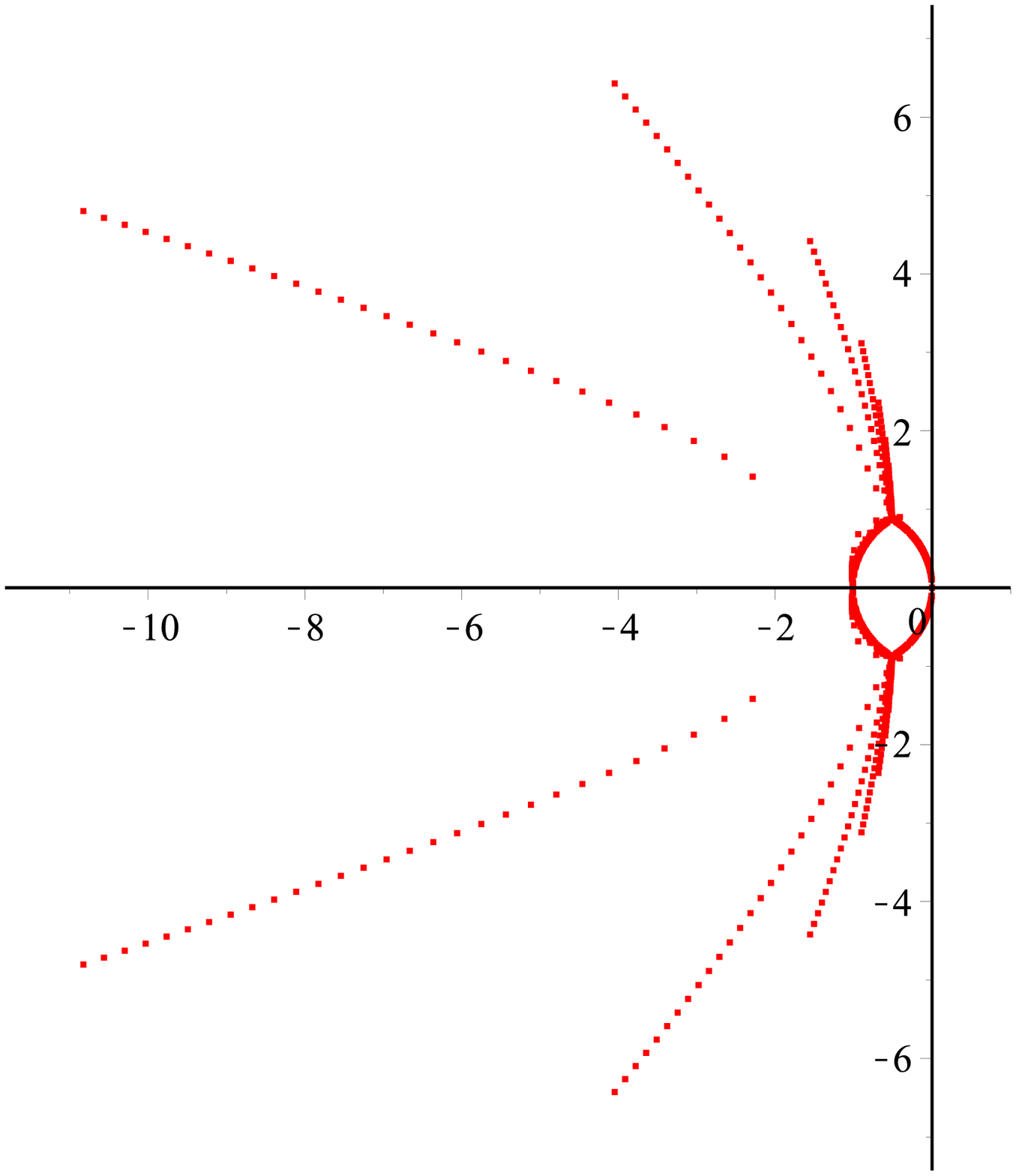}
	\includegraphics[height=6.4cm]{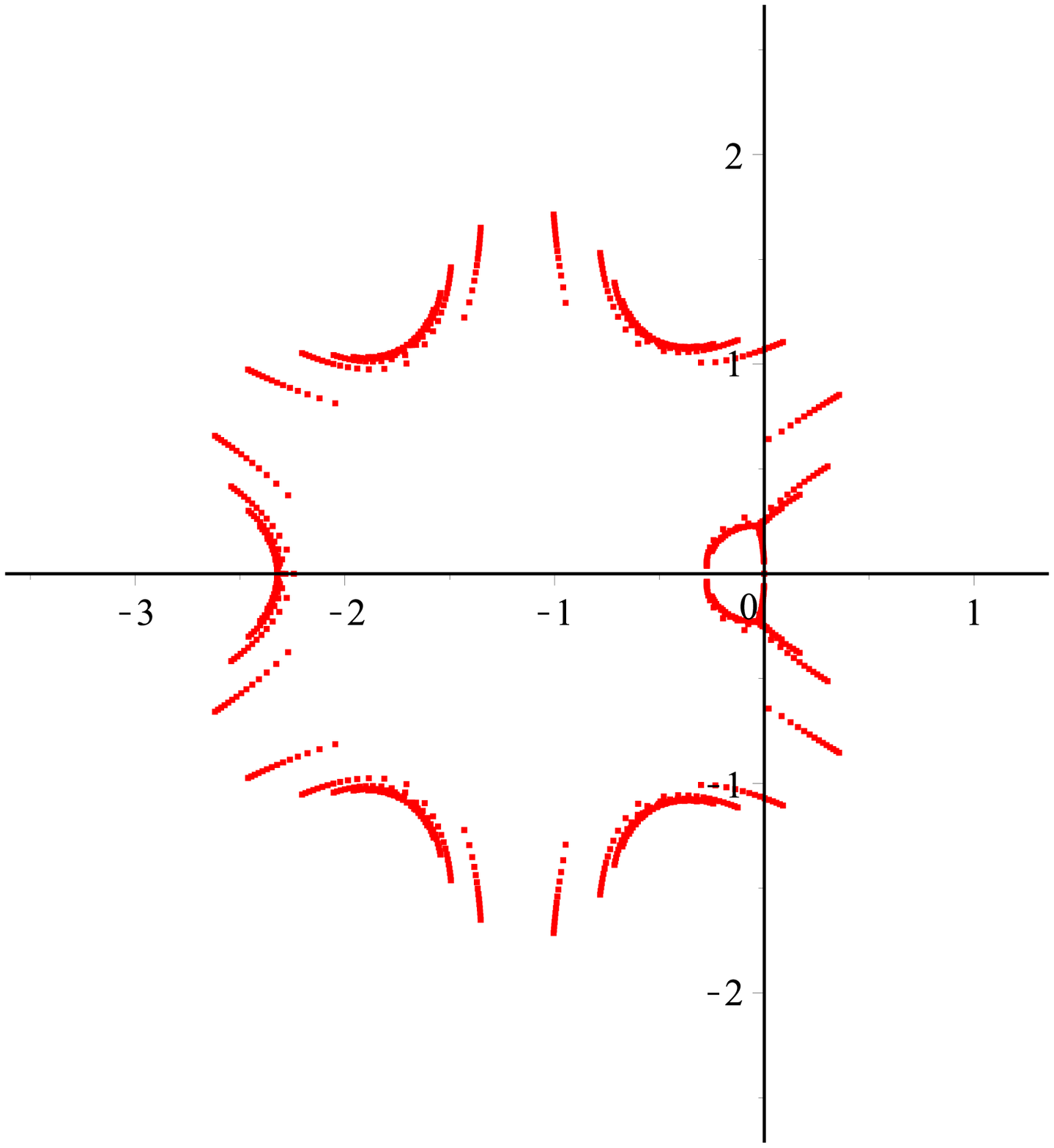}
	\caption{Total domination roots of $K_{1,n}[K_2]$ and $K_{1,n}[K_7]$, for $2 \leq n \leq 30$, respectively.}\label{K_1,n[K_m]}
\end{figure}

\begin{theorem}\label{join}
	Let $G$ and $H$ be two graphs of order $m$ and $n$, respectively. The total domination polynomial of join of these two graphs is
	\[
	D_t(G\vee H)=((1+x)^m-1)((1+x)^n-1)+D_t(G,x)+D_t(H,x).
	\]
\end{theorem}

\begin{theorem}
	For every natural numbers $m,n$, the total domination polynomial of $K_{1,n}[K_m]$ is
	\[
	D_t(K_{1,n}[K_m],x)=(1+x)^{mn}((1+x)^m-1)+((1+x)^m-mx-1)^n-mx.
	\]
\end{theorem}
\proof
For two natural numbers $m,n$, $K_{1,n}[K_m]\cong k_m \vee nK_m$. So by Theorem \ref{join}, it is easy to see the equation is true.\qed

Using Maple we think that for two natural numbers $m$,$n$, if $m$ and $n$ are even or $n$ is odd, then the total domination polynomial of $K_{1,n}[K_m]$ has no real roots. However, until now all attempts to prove this failed.
See the total domination roots of $K_{1,n}[K_2]$ and $K_{1,n}[K_7]$ for $2 \leq n \leq 30$ in Figure \ref{K_1,n[K_m]}.
        



\begin{thebibliography}{1}
	
	\bibitem{euro}{\sc S.\ Akbari {\rm and} S.\ Alikhani {\rm and} Y.\ H.\ Peng}: \textit{ Characterization of graphs using domination polynomials}. Eur. J. Combin. {\bf 31} (2010), 1714--1724. 
	
	
	
	
	
	
	
	\bibitem{nasrin3}{\sc S.\ Alikhani {\rm and} N.\ Jafari}: \textit{ Some new results on the total domination polynomial of a graph}. Ars Combin. In press. Available at \texttt{http://arxiv.org/abs/1705.00826}.
	
	
	
	\bibitem{bridges}{\sc W.\ G.\ Bridges {\rm and} R.\ A.\ Mena}: \textit{ Multiplicative cones- a family of three eigenvalue graph}. Aequationes Math. {\bf 22} (1981), 208--214.
	
	
	\bibitem{Brown} {\sc J.\ I.\ Brown {\rm and} J.\ Tufts}:\textit{ On the roots of domination polynomials}. Graphs Combin. {\bf 30} (2014), 527--547.
	
	\bibitem{BHN} {\sc J.\ I.\ Brown {\rm and} C.\ A.\ Hickman{\rm and} R.\ J.\ Nowakowski}: \textit{ On the location of roots of independence polynomials}. J. Algebraic Combin. {\bf 19} (2004), 273--282. 
	
	\bibitem{dam1} {\sc E.\ R.\ Van Dam}: \textit{ Regular graphs with four eigenvalues}.Linear Algebra Appl, {\bf 226/228} (1995), 139--162.
	
	\bibitem{dam2} {\sc E.\ R.\ Van Dam}: \textit{ Graphs with Few Eigenvalues}, An Interplay between Combinatorics and Algebra, Center Dissertation Series 20, Thesis, Tilburg University, 1996.
	
	\bibitem{dam3} {\sc E.\ R.\ Van Dam}: \textit{ Nonregular graphs with three eigenvalues}. J. Combin. Theory Ser, B {\bf 73} (1998), 101--118.
	
	\bibitem{dam4} {\sc E.\ R.\ Van Dam {\rm and} W.\ H.\ Haemers}: \textit{ Which graphs are determined by their spectrum?}. Linear Algebra Appl, {\bf 373} (2003), 241--272.
	
	\bibitem{Dod} {\sc M.\ Dod}: \textit{ The total domination polynomial and its generalization}. In: Congressus Numerantium, {\bf 219} (2014), 207--226.
	
	
	\bibitem{15} {\sc C.\ D.\ Godsil}: \textit{ Algebraic Combinatorics}. Chapmanand Hall, NewYork. 1993.
	
	\bibitem{Harary} {\sc F.\ Harary}: \textit{ On the group of the composition of two graphs}. Duke Math. J.{\bf 26} (1959), 29--36.
	
	\bibitem{Harar} {\sc F.\ Harary}: \textit{ Graph Theory}. Addison-Wesley, Reading, MA (1969).
	
	
	\bibitem{Heil} {\sc O.\ J.\ Heilmann {\rm and} E.\ H.\ Lieb}: \textit{ Theory of monomer-dimer systems}, Comm. Math. Phys. {\bf 25} (1972), 190--232. 
	
	\bibitem{Hening}{\sc M.\ A.\ Henning {\rm and} A.\ Yeo}: \textit{Total domination in graphs }. Springer Monographs in
	Mathematics, 2013.
	
	
	\bibitem{nasrin}{\sc N.\ Jafari  {\rm and} S. Alikhani}:\textit{ On the roots of total domination polynomial of graphs}, J. Discrete Math. Sci. Crypt., to appear,  Available at \texttt{http://arxiv.org/abs/1605.02222}.
	
	\bibitem{muz} {\sc M.\ Klin {\rm and} M.\ Muzychuk}: \textit{ On graphs with three eigenvalues}.
	Discrete Math. {\bf 189} (1998), 191--207.
	
	
	\bibitem{omidi} {\sc H.\ Chuang {\rm and} G.\ R.\ Omidi}:\textit{ Graphs with three distinct eigenvalues and largest
		eigenvalue less than 8}. Linear Algebra Appl. {\bf 430} (2009), 2053--2062.
	
	\bibitem{Schiermeyer} {\sc Z.\ Ryj\'{a}\v{c}ek {\rm and} I.\ Schiermeyer}: \textit{ The flower conjecture in special classes of graphs}. Discuss. Math. Graph Theory, {\bf 15} (1995), 179--184.
	
	
	\bibitem{Oboudi}{\sc M.\ R.\ Oboudi}: \textit{ On the roots of domination polynomial of graphs}. Discrete Appl. Math. {\bf 205} (2016), 126--131. 
	
	
	
	
\end{thebibliography}
\end{document}